\documentclass[11pt,fleqn,twoside]{article}
\usepackage{amsfonts,amssymb,latexsym}
\makeatletter
\newcommand{\prava}{\footnotesize\it
\begin{flushright}
\begin{minipage}{18cm}
Copyright \copyright 1998 by P. Wall
\end{minipage}
\end{flushright}}

\newcommand{\name}[1]{\begin{flushleft}
                       \LARGE \bf #1
                       \end{flushleft}\vspace{-3mm}}

\newcommand{\Author}[1]{\begin{flushleft}
                       \it #1 \end{flushleft}}

\newcommand{\Adress}[1]{\begin{flushleft}
                       \it #1 \end{flushleft}}

\newcommand{\Date}[1]{\begin{flushleft}
                      \small  \it #1 \end{flushleft}}

\newcommand{\ehkol}{Author \ name}
\newcommand{\ohkol}{Article \ name}
\renewcommand{\@evenhead}{
\hspace*{-3pt}\raisebox{-15pt}[\headheight][0pt]{\vbox{\hbox to \textwidth
{\thepage \hfil \ehkol}\vskip4pt \hrule}}}
\renewcommand{\@oddhead}{
\hspace*{-3pt}\raisebox{-15pt}[\headheight][0pt]{\vbox{\hbox to \textwidth
{\ohkol \hfil \thepage}\vskip4pt\hrule}}}
\renewcommand{\@evenfoot}{}
\renewcommand{\@oddfoot}{}

\newcommand{\be}{\begin{equation}}
\newcommand{\ee}{\end{equation}}
\newcommand{\ba}{\hspace*{-5pt}\begin{array}}
\newcommand{\ea}{\end{array}}

\newcommand{\ds}{\displaystyle}
\makeatother

\begin{document}

\newcommand{\limfunc}[1]{\mathop{\rm #1}\nolimits}
\setcounter{page}{331}
\thispagestyle{empty}

\renewcommand{\ehkol}{P. Wall}
\renewcommand{\ohkol}{Homogenization 
for Nonlinear Monotone Operators}

\begin{flushleft}
\footnotesize \sf
Journal of Nonlinear Mathematical Physics \qquad 1998, V.5, N~3,
\pageref{wall-fp}--\pageref{wall-lp}.
\hfill {\sc Article}
\end{flushleft}

\vspace{-5mm}

\renewcommand{\footnoterule}{}
{\renewcommand{\thefootnote}{}
 \footnote{\prava}}

\name{Some Homogenization and Corrector Results\\
for Nonlinear Monotone Operators}\label{wall-fp}

\Author{Peter WALL}

\Adress{Department of Mathematics, Lule\aa\
University of Technology, SE-971 87 Lule\aa, Sweden}

\Date{Received March 5, 1998; Accepted May 9, 1998}

\begin{abstract}
\noindent
This paper deals with the limit behaviour of the solutions of quasi-linear
equations of the form \ $\ds -\limfunc{div}\left(a\left(x,
x/{\varepsilon _h},Du_h\right)\right)=f_h$
on $\Omega $ with Dirichlet boundary conditions. The sequence $(\varepsilon
_h)$ tends to $0$ and the map $a(x,y,\xi )$ is periodic in $y$,
monotone in $\xi $ and satisf\/ies suitable continuity conditions. It
is proved that $u_h\rightarrow u$ weakly in $H_0^{1,2}(\Omega )$,
where $u$ is the solution of a homogenized problem \
$-\limfunc{div}(b(x,Du))=f$ on $\Omega $. We also prove some
corrector results, i.e. we f\/ind $(P_h)$ such that
$Du_h-P_h(Du)\rightarrow 0$ in $L^2(\Omega ,R^n)$.
\end{abstract}

\section{Introduction}
In mathematical models of microscopically non-homogeneous media various
local characteristics are usually described by functions of the form
$\ds a\left(x/{\varepsilon _h}\right)$ where $\varepsilon _h>0$
is a small parameter.
The function $a(x)$ can be periodic or belong to some other
specif\/ic class. To compute the properties of a micro
non-homogeneous medium is an extremely dif\/f\/icult task since the
coef\/f\/icients are rapidly oscillating functions.
Therefore, it is necessary to apply asymptotic analysis to the problems of
micro non-homogeneous media, which immediately leads to the concept of
homogenization. When the parameter $\varepsilon _h$ is very small the medium
will act as a homogeneous medium. To characterize this homogeneous medium is
one of the main tasks in the homogenization theory. For more information
concerning the homogenization theory the reader is referred to
\cite{bensoussan1978}, \cite{jikov1994a} and \cite{perssonetal1993}.
In this paper we consider the homogenization problem for monotone
operators and the local behavior of the solutions. Monotone operators
are very important in the study of nonlinear partial dif\/ferential
equations. The problem we study here can be used to model dif\/ferent
nonlinear stationary conservation laws, e.g. stationary temperature
distribution. For a more detailed discussion concerning dif\/ferent
applications see \cite{zeidler19904}.

We will study the limit behavior of the sequence of solutions $(u_h)$
to the Dirichlet boundary value problem
\[
\left\{
\begin{array}{l}
\ds -\limfunc{div}\left(a\left(x,\frac x{\varepsilon _h},Du_h\right
)\right )=f_h\quad \mbox{on}
\quad \Omega , \\[5mm]
u_h\in H_0^{1,2}(\Omega ),
\end{array}
\right.
\]
where $f_h\rightarrow f$ in $H^{-1,2}(\Omega )$ and $\varepsilon
_h\rightarrow 0$. Moreover, the map $a(x,y,\xi )$ is def\/ined on $\Omega
\times R^n\times R^n$ and is assumed to be periodic in $y$, uniformly
Lipschitz continuous in $\xi $ and uniformly monotone in $\xi $. We also
need some continuity restriction on $a(\cdot ,y,\xi )$. We will consider two
dif\/ferent cases, namely when $a(x,y,\xi )$ is of the form
\[
a(x,y,\xi )=\sum_{i=1}^N\chi _{\Omega _i}(x)a_i(y,\xi ),
\]
or when $a(x,y,\xi )$ satisf\/ies that
\[
\left| a(x_1,y,\xi )-a(x_2,y,\xi )\right| ^2\leq \omega (\left|
x_1-x_2\right| )\left| \xi \right| ^2,
\]
where $\omega :R\rightarrow R$ is continuous, increasing and $\omega (0)=0.$
In both cases we will prove that $u_h\rightarrow u$ weakly in
$H^{1,2}(\Omega )$ and that $u$ is the solution of the
{\it homogenized}  problem
\[
\left\{
\begin{array}{l}
-\limfunc{div}(b(x,Du))=f\quad \mbox{on} \quad \Omega , \\[2mm]
u\in H_0^{1,2}(\Omega ).
\end{array}
\right.
\]
We will prove that the operator $b$ has the same structure properties as $a$
and is given by
\[
b(x,\xi )=\int_Ya\left(x,y,\xi +Dv^{\xi ,x}(y)\right)\,dy,
\]
where $v^{\xi ,x}$ is the solution of the {\it cell-problem}
\begin{equation}
\left\{
\begin{array}{l}
-\limfunc{div}\left(a\left(x,y,\xi +Dv^{\xi ,x}(y)\right)\right)=0
\quad \mbox{on}
\quad Y, \\[3mm]
v^{\xi ,x}\in H_{\square }^{1,2}(Y),
\end{array}
\right.  \label{cellproblem0}
\end{equation}
where $Y$ is a cell of periodicity and $H_{\square }^{1,2}(Y)$ is the subset
of $H^{1,2}(Y)$ such that $u$ has mean value $0$ and $u$ is $Y$-periodic.
The homogenization problem for monotone operators of this type has been
studied by several authors but with no dependence in $x$, i.e. $a$ is on the
form $a(x,y,\xi )=a(y,\xi )$. Here we mention \cite{fusco1986} where the
problem was studied in the Sobolev space $H^{1,p}$, $1<p<\infty $, with
appropriate continuity and monotonicity conditions. In \cite{chiadopiat1990}
the corresponding multi-valued case was considered. Moreover, the almost
periodic case was treated in \cite{braides1992}.

The weak convergence of $u_h$ to $u$ in $H^{1,2}(\Omega )$ implies that $%
u_h-u\rightarrow 0$ in $L^2(\Omega )$ but in general we only have that $%
Du_h-Du\rightarrow 0$ weakly in $L^2(\Omega ,R^n)$. However, we will prove
that it is possible to express $Du_h$ in terms of $Du$, up to a rest which
converges strongly in $L^2(\Omega ,R^n)$. This is done by constructing a
family of correctors $P_h(x,\xi ,t)$, def\/ined by
\begin{equation}
P_h(x,\xi ,t)=P\left(\frac x{\varepsilon _h},\xi ,t\right)=\xi
+Dv^{\xi ,t} \left(\frac x{\varepsilon _h}\right).  \label{pdef}
\end{equation}

Let $(M_h)$ be a family of linear operators converging to the identity map
on $L^2(\Omega ,R^n)$ such that $M_hf$ is a step function for every $f\in
L^2(\Omega ,R^n)$. Moreover, let $\gamma _h$ be a step function
approximating the identity map on $\Omega $. We will show that
\[
Du_h-P_h(x,M_hDu,\gamma _h)\rightarrow 0\quad \mbox{in} \quad
L^2(\Omega ,R^n).
\]

Results concerning correctors for an even more general case than studied in
\cite{fusco1986} can be found in \cite{dalmaso1992}. The almost periodic
case has been studied in \cite{braides1991}.

\section{Preliminaries and notations}

Let $\Omega $ be a open bounded subset of $R^n$, $\left| E\right| $ denotes
the Lebesgue measure in $R^n$ and $(\cdot ,\cdot )$ is the Euclidean scalar
product on $R^n$. Let $\left\{ \Omega _i\subset \Omega :i=1,\ldots
,N\right\} $ be a family of disjoint open sets such that $\left| \Omega
\backslash \cup _{i=1}^N\Omega _i\right| =0$ and $\left| \partial \Omega
_i\right| =0$. Let $\left( \varepsilon _h\right) $ be a decreasing sequence
of real numbers such that $\varepsilon _h\rightarrow 0$ as $h\rightarrow
\infty $. $Y=\left( 0,1\right) ^n$ is the unit cube in $R^n$ and
$Y_h^j=\varepsilon _h(j+Y)$, where $j\in Z^n$, i.e. the translated
image of $\varepsilon _hY$ by the vector $\varepsilon _hj$. We also
def\/ine the following index sets:
\[
\ba{l}
\ds J_h =\left\{ j\in Z^n:\overline{Y}_h^j\subset \Omega \right\}
,\;J_h^i=\left\{ j\in Z^n:\overline{Y}_h^j\subset \Omega _i\right\} , \\[3mm]
\ds B_h^i =\left\{ j\in Z^n:\overline{Y}_h^j\cap \Omega _i\neq \emptyset ,\;
\overline{Y}_h^j\backslash \Omega _i\neq \emptyset \right\} .
\ea
\]
Moreover, we def\/ine $\Omega _i^h=\cup _{j\in J_h^i}\overline{Y}_h^j$ and $%
F_i^h=\cup _{j\in B_h^i}Y_h^j$.

In a corresponding way let $\left\{ \Omega _i^k\subset \Omega :i\in
I_k\right\} $ denote a family of disjoint open sets with diameter
less than $\ds \frac 1k$ such that $\left| \Omega \backslash \cup _{i\in
I_k}\Omega
_i^k\right| =0$ and $\left| \partial \Omega _i^k\right| =0$. We also def\/ine
the following index sets:
\[
\ba{l}
\ds J_h^{i,k} =\left\{ j\in Z^n:\overline{Y}_h^j\subset \Omega _i^k\right\} ,
\\[3mm]
\ds B_h^{i,k} =\left\{ j\in Z^n:\overline{Y}_h^j\cap \Omega _i^k\neq
\emptyset ,\;\overline{Y}_h^j\backslash \Omega _i^k\neq \emptyset \right\} .
\ea
\]
Let $\Omega _i^{k,h}=\cup _{j\in J_h^{i,k}}\overline{Y}_h^j$ and
$F_i^{k,h}=\cup _{j\in B_h^{i,k}}Y_h^j$.

Corresponding to $f\in L^2(\Omega ,R^n)$ we def\/ine the function $%
M_hf:R^n\rightarrow R^n$ by
\[
(M_hf)(x)=\sum_{j\in J_h}\chi _{Y_h^j}(x)\xi _h^j,
\]
where $\ds \xi _h^j=\frac 1{\left| Y_h^j\right| }\int_{Y_h^j}f\,dx$ and
$\chi _E$ is the characteristic function of the set $E$ (in order to
def\/ine $\xi _h^j$
for all $j\in Z^n$ we treat $f$ as $f=0$ outside $\Omega $). It is well
known that
\begin{equation}
M_hf\rightarrow f\quad \mbox{in}\quad L^2(\Omega ,R^n),  \label{approx}
\end{equation}
see \cite{royden1988}, page 129. We also def\/ine the step function $\gamma
_h:\Omega \rightarrow \Omega $ by
\begin{equation}
\gamma _h(x)=\sum_{j\in J_h}\chi _{Y_h^j}(x)x_h^j,  \label{gamma}
\end{equation}
where $x_h^j\in Y_h^j.$ Moreover, $C$ will be a constant that may dif\/fer
from one place to an other. Let $a:\Omega \times R^n\times R^n\rightarrow
R^n $ be a function such that $a(x,\cdot ,\xi )$ is Lebesgue
measurable and $Y$-periodic for $x\in \Omega $ and $\xi \in R^n.$ We
also assume that $a$
satisf\/ies the following monotonicity and continuity conditions: There exists
two constants $0<\alpha \leq \beta <\infty $ such that
\be
(a(x,y,\xi _1)-a(x,y,\xi _2),\xi _1-\xi _2) \geq \alpha \left| \xi _1-\xi
_2\right| ^2,  \label{amon}
\ee
\be
\left| a(x,y,\xi _1-a(x,y,\xi _2)\right| \leq \beta \left| \xi _1-\xi
_2\right| ,  \label{acont}
\ee
for $x\in \Omega,$ a.e. $y\in R^n$ and every $\xi \in R^n$. Moreover we
assume that
\begin{equation}
a(x,y,0)=0,  \label{axy0}
\end{equation}
for $x\in \Omega ,$ a.e. $y\in R^n$. Let $(f_h)$ be a sequence in
$H^{-1,2}(\Omega )$ which converges to $f$.

The solution $v^{\xi ,x}$ of the cell-problem (\ref{cellproblem0}) can be
extended by periodicity to an element in $H_{\mbox{\scriptsize {\rm
loc}}}^{1,2}(R^n)$, still denoted by $v^{\xi ,x}$ and
\begin{equation}
\int_{R^n}\left(a\left(x,y,\xi +Dv^{\xi ,x}(y)\right),
D\phi (y)\right)\,dy=0\quad \mbox{for
every} \quad \phi \in C_0^\infty (R^n).  \label{divwhk}
\end{equation}
The following compensated compactness lemma will be used frequently, see
\cite{jikov1994a}, page 4.

\medskip

\noindent
{\bf Lemma 1.}\label{compensatedcompctness}
{\it Let $1<p<\infty $. Moreover, let $(v_h)$ be a
sequence in $L^q(\Omega ,R^n)$ which converges weakly to $v$,
$(-\limfunc{div}v_h)$ converges to $-\limfunc{div}v$ in
$W^{-1,q}(\Omega )$ and let $(u_h)$
be a sequence which converges weakly to $u$ in $W^{1,p}(\Omega ).$ Then
\[
\int_\Omega (v_h,Du_h)\phi \,dx\rightarrow \int_\Omega (v,Du)\phi \,dx,
\]
for every $\phi \in C_0^\infty (\Omega ).$
}

\section{Some homogenization results}

Let $a(x,y,\xi )$ satisfy (\ref{amon}), (\ref{acont}), (\ref{axy0}) and one
of the following conditions:

\begin{enumerate}
\item[(i)]  $a$ is on the form
\begin{equation}
a(x,y,\xi )=\sum_{i=1}^N\chi _{\Omega _i}(x)a_i(y,\xi ).  \label{astep}
\end{equation}

\item[(ii)]  there exist a function $\omega :R\rightarrow R$ which is
continuous, increasing and $\omega (0)=0$ such that
\begin{equation}
\left| a(x_1,y,\xi -a(x_2,y,\xi )\right| ^2\leq \omega (\left|
x_1-x_2\right| )\left| \xi \right| ^2,  \label{cont3}
\end{equation}
for $x\in \Omega ,$ a.e. $y\in R^n$ and every $\xi \in R^n$.
\end{enumerate}

Now we consider the weak Dirichlet boundary value problems, one for each $h,$
\begin{equation}
\left\{
\begin{array}{l}
\ds \int_\Omega \left(a\left(x,\frac x{\varepsilon _h},Du_h
\right),D\phi \right)\,dx=\left\langle
f_h,\phi \right\rangle \quad \mbox{for every } \quad \phi \in
H_0^{1,2}(\Omega ), \\[5mm]
u_h\in H_0^{1,2}(\Omega ).
\end{array}
\right.  \label{weakproblem}
\end{equation}
By a standard result in the existence theory for boundary value problems
def\/ined by monotone operators these problems have unique solution for
each $h $, see e.g. \cite{zeidler19902b}. Furthermore, by choosing
$\phi =u_h$ in (\ref{weakproblem}), taking into account (\ref{amon}),
(\ref{axy0}) and using the fact that $(f_h)$ is bounded we have that
\[
\ba{l}
\ds \alpha \int_\Omega \left| Du_h\right| ^2\,dx \leq \int_\Omega
\left(a\left(x,\frac
x{\varepsilon _h},Du_h\right),Du_h\right)\,dx=\left\langle
f_h,u_h\right\rangle \\[5mm]
\phantom{\alpha \int_\Omega \left| Du_h\right| ^2\,dx}
\leq \left\| f_h\right\| _{H^{-1,2}(\Omega )}\left\| u_h\right\|
_{H_0^{1,2}(\Omega )}\leq C\left\| u_h\right\| _{H_0^{1,2}(\Omega )},
\ea
\]
where $C$ does not depend on $h$. The Poincar\'{e} inequality then implies
that
\begin{equation}
\left\| u_h\right\| _{H_0^{1,2}(\Omega )}\leq C,  \label{uhbounded}
\end{equation}
where $C$ does not depend on $h$. Therefore there exists a
subsequence $(h^{\prime })$ such that
\begin{equation}
u_{h^{\prime }}\rightarrow u_{*}\quad \mbox{weakly in} \quad H_0^{1,2}(\Omega ).
\label{uhconv}
\end{equation}
It is now natural to raise the following question: does $u_{*}$ satisfy an
equation of the same type as that satisf\/ied by $u_h$? The answer to this
question is given in the following theorems:

\medskip

\noindent
{\bf Theorem 1.} \label{transth}
{\it Let $a$ satisfy (\ref{amon}), (\ref{acont}), (\ref{axy0}) and
(\ref{astep}). Moreover, let $(u_h)$ be the solutions
of~(\ref{weakproblem}). Then
\[
u_h\rightarrow u\quad \mbox{weakly in} \quad H_0^{1,2}(\Omega ),
\]
\begin{equation}
a\left(x,\frac x{\varepsilon _h},Du_h\right)\rightarrow b(x,Du) \quad
\mbox{weakly in} \quad
L^2(\Omega ;R^n),  \label{ahbounded}
\end{equation}
where $u$ is the unique solution of the homogenized problem
\begin{equation}
\left\{
\begin{array}{l}
\ds\int_\Omega (b(x,Du),D\phi )\,dx=\left\langle f,\phi \right\rangle
\quad \mbox{for every}\quad \phi \in H_0^{1,2}(\Omega ), \\[5mm]
u\in H_0^{1,2}(\Omega ).
\end{array}
\right.   \label{homogenizedproblem}
\end{equation}
The operator $b:\Omega \times R^n\rightarrow R^n$ is defined a.e. as
\[
b(x,\xi )=\sum_{i=1}^N\chi _{\Omega _i}(x)\int_Ya_i\left(y,\xi +Dv^{\xi
,x_i}(y)\right)\,dy=\sum_{i=1}^N\chi _{\Omega _i}(x)b_i(\xi ),
\]
where $x_i\in \Omega _i$ and $v^{\xi ,x_i}$ is the unique solution of the
cell problem
\begin{equation}
\left\{
\begin{array}{l}
\ds \int_Y\left(a_i\left(y,\xi +Dv^{\xi ,x_i}(y)\right),D\phi (y)\right)\,dy=0
\quad \mbox{for every}\quad \phi
\in H_{\square }^{1,2}(Y), \\[5mm]
v^{\xi ,x_i}\in H_{\square }^{1,2}(Y).
\end{array}
\right.   \label{cellproblem}
\end{equation}
}

\noindent
{\bf Proof.}
The proof follows by using the ideas in \cite{defranceschi1993} where the
case $N=1$ is treated for the details the reader is referred to
\cite{wall1998c}.  \hfill {}$\square $

\medskip

\noindent
{\bf Theorem 2.} \label{cth}
{\it Let $a$ satisfy (\ref{amon}), (\ref{acont}), (\ref{axy0}) and
(\ref{cont3}). Moreover, let $(u_h)$ be the solutions
of~(\ref{weakproblem}).
Then
\[
u_h \rightarrow u \quad \mbox{weakly in} \quad H_0^{1,2}(\Omega ),
\]
\[
a\left(x,\frac x{\varepsilon _h},Du_h\right) \rightarrow b(x,Du)\quad
\mbox{weakly in} \quad L^2(\Omega ;R^n),
\]
where $u$ is the unique solution of
\[
\left\{
\begin{array}{l}
\ds \int_\Omega (b(x,Du),D\phi )\,dx=\left\langle f_h,\phi \right\rangle
\quad \mbox{for every} \quad \phi \in H_0^{1,2}(\Omega ), \\[5mm]
u\in H_0^{1,2}(\Omega ).
\end{array}
\right.
\]
The operator $b:\Omega \times R^n\rightarrow R^n$ is defined as
\[
b(x,\xi )=\int_Ya\left(x,y,\xi +Dv^{\xi ,x}(y)\right)\,dy,
\]
where $v^{\xi ,x}$ is the unique solution of the cell-problem
\begin{equation}
\left\{
\begin{array}{l}
\ds \int_Y\left(a\left(x,y,\xi +Dv^{\xi ,x}(y)\right),D\phi
\right)\,dy=0 \quad \mbox{for
every} \quad \phi \in
H_{\square }^{1,2}(Y), \\[5mm]
v^{\xi ,x}\in H_{\square }^{1,2}(Y).
\end{array}
\right.   \label{cellproblem2}
\end{equation}
}

Before we prove this theorem we make some def\/initions and prove some lemmas
that will be useful in the proof. Def\/ine the function
\[
a^k(x,y,\xi ):=\sum_{i\in I_k}\chi _{\Omega _i^k}(x)a(x_i^k,y,\xi ),
\]
where $x_i^k\in \Omega _i^k$. Consider the boundary value problems
\begin{equation}
\left\{
\begin{array}{l}
\ds \int_\Omega \left(a^k
\left(x,\frac x{\varepsilon _h},Du_h^k\right),D\phi \right)\,dx=
\left\langle
f_h,\phi \right\rangle \quad \mbox{for every}\quad \phi \in
H_0^{1,2}(\Omega ), \\[5mm]
u_h^k\in H_0^{1,2}(\Omega ).
\end{array}
\right.  \label{uhk}
\end{equation}
The conditions for Theorem~1 are satisf\/ied and the theorem
implies that there exists a $u_{*}^k$ such that
\[
u_h^k\rightarrow u_{*}^k \quad \mbox{weakly in} \quad
H_0^{1,2}(\Omega ) \quad \mbox{as} \quad
h\rightarrow \infty ,
\]
and $u_{*}^k$ is the unique solution of
\begin{equation}
\left\{
\begin{array}{l}
\ds \int_\Omega \left(b^k\left(x,Du_{*}^k\right),D\phi
\right)\,dx=\left\langle f,\phi \right\rangle
\quad \mbox{for every} \quad \phi \in H_0^{1,2}(\Omega ), \\[5mm]
u_{*}^k\in H_0^{1,2}(\Omega ),
\end{array}
\right.  \label{ustark}
\end{equation}
where
\[
b^k(x,\xi )=\sum_{i\in I_k}\chi _{\Omega _i^k}(x)\int_Ya\left(x_i^k,y,\xi
+Dv^{\xi ,x_i^k}(y)\right)\,dy=\sum_{i\in I_k}\chi _{\Omega
_i^k}(x)b\left(x_i^k,\xi \right),
\]
where $v^{\xi ,x_i^k}$ is the solution of
\[
\left\{
\begin{array}{l}
\ds \int_Y\left(a\left(x_i^k,y,\xi +Dv^{\xi ,x_i^k}(y)\right),D\phi
\right) \,dy=0
\quad \mbox{for every} \quad
\phi \in H_{\square }^{1,2}(Y), \\[5mm]
v^{\xi ,x_i^k}\in H_{\square }^{1,2}(Y).
\end{array}
\right.
\]

\noindent
{\bf Proof of Theorem~2.} First we prove that
$u_h\rightarrow u$ weakly
in $H_0^{1,2}(\Omega )$. If $g\in H^{-1,2}(\Omega ),$ then
\[\hspace*{-9pt}
\ba{l}
\ds \lim_{h\rightarrow \infty }\left\langle g,u_h-u\right\rangle
=\lim_{k\rightarrow \infty }\lim_{h\rightarrow \infty }\left\langle
g,u_h-u\right\rangle =\lim_{k\rightarrow \infty }\lim_{h\rightarrow \infty
}\left\langle g,u_h-u_h^k+u_h^k-u_{*}^k+u_{*}^k-u\right\rangle \\[4mm]
\ds \qquad \leq \lim_{k\rightarrow \infty }\lim_{h\rightarrow \infty }\left\|
g\right\| _{H^{-1,2}(\Omega )}\left\| u_h-u_h^k\right\| _{H_0^{1,2}(\Omega
)}+\lim_{k\rightarrow \infty }\lim_{h\rightarrow \infty }\left\langle
g,u_h^k-u_{*}^k\right\rangle \\[4mm]
\qquad +\lim_{k\rightarrow \infty }\left\| g\right\|
_{H^{-1,2}(\Omega )}\left\| u_{*}^k-u\right\| _{H_0^{1,2}(\Omega )}.
\ea
\]
It is enough to prove that all three terms on the right hand side are zero.
We do this in three steps.

\textbf{Step 1.} Let us prove that
\begin{equation}
\lim_{k\rightarrow \infty }\lim_{h\rightarrow \infty }\left\|
u_h-u_h^k\right\| _{H_0^{1,2}(\Omega )}=0  \label{step1a}
\end{equation}
By def\/inition
\[
\int_\Omega \left(a^k\left(x,\frac x{\varepsilon
_h},Du_h^k\right),D\phi \right)\,dx
=\left\langle f_h,\phi \right\rangle \quad \mbox{for every} \quad \phi \in
H_0^{1,2}(\Omega ),
\]
\[
\int_\Omega \left(a\left(x,\frac x{\varepsilon _h},Du_h\right),D\phi
\right) \,dx =\left\langle
f_h,\phi \right\rangle \quad \mbox{for every} \quad \phi \in
H_0^{1,2}(\Omega ).
\]
This implies that we for $\phi =u_h^k-u_h$ have
\[
\ba{l}
\ds \int_\Omega \left(a^k\left(x,\frac x{\varepsilon
_h},Du_h^k\right) -a^k\left(x,\frac
x{\varepsilon _h},Du_h\right),Du_h^k-Du_h\right)\,dx \\[5mm]
\ds \qquad =\int_\Omega \left(a\left(x,\frac x{\varepsilon
_h},Du_h\right)- a^k\left(x,\frac x{\varepsilon
_h},Du_h\right),Du_h^k-Du_h\right)\,dx.
\ea
\]
According to the monotonicity of $a,$ (\ref{amon}), the Schwarz inequality
and the H\"{o}lder inequality we obtain that
\[
\ba{l}
\ds \alpha \int_\Omega \left| Du_h^k-Du_h\right| ^2\,dx \\[5mm]
\ds \qquad \leq \left( \int_\Omega \left| a\left(x,\frac
x{\varepsilon _h},Du_h\right)-a^k\left(x,\frac x{\varepsilon
_h},Du_h\right) \right| ^2\,dx\right) ^{\frac
12}\left( \int_\Omega \left| Du_h^k-Du_h\right| ^2\,dx\right) ^{\frac 12}
\ea
\]
i.e.
\[
\left\| Du_h^k-Du_h\right\| _{L^2(\Omega ;R^n)}^2\leq \frac 1{\alpha
^2}\int_\Omega \left| a\left(x,\frac x{\varepsilon
_h},Du_h\right)-a^k \left(x,\frac x{\varepsilon _h},Du_h\right)\right| ^2\,dx.
\]
Thus, in view of the continuity condition (\ref{cont3}) on $a$, we get
\begin{equation}
\left\| Du_h^k-Du_h\right\| _{L^2(\Omega ;R^n)}^2\leq \frac 1{\alpha
^2}\omega \left(\frac 1k\right)\int_\Omega \left| Du_h\right|
^2\,dx\leq \frac C{\alpha
^2}\omega \left(\frac 1k\right),  \label{bound}
\end{equation}
where we in the last inequality used the fact that there exists a constant $%
C $ independent of $h$ such that $\left\| Du_h\right\| _{L^2(\Omega
;R^n)}^2\leq C$. Since $\left\| D\cdot \right\| _{L^2(\Omega ;R^n)}^2$ is an
equivalent norm on $H_0^{1,2}(\Omega )$~(\ref{bound}) implies that
\[
\left\| u_h^k-u_h\right\| _{H_0^{1,2}(\Omega )}\rightarrow 0,
\]
as $k\rightarrow \infty $ uniformly in $h$. This means that we can change
the order in the limit process in (\ref{step1a}) and (\ref{step1a}) follows
by taking (\ref{bound}) into account.
\hfill {}$\square $

\textbf{Step 2.} We observe that
\[
\lim_{k\rightarrow \infty }\lim_{h\rightarrow \infty }\left\langle
g,u_h^k-u_{*}^k\right\rangle =0,
\]
as a direct consequence of Theorem~1. \hfill {}$\square $

\textbf{Step 3.} Let us prove that
\begin{equation}
\lim_{k\rightarrow \infty }\left\| u_{*}^k-u\right\| _{H_0^{1,2}(\Omega )}=0.
\label{step3}
\end{equation}
By def\/inition we have that
\[
\int_\Omega \left(b^k\left(x,Du_{*}^k\right),D\phi \right)\,dx =
\left\langle f,\phi
\right\rangle \quad \mbox{for every} \quad \phi \in
H_0^{1,2}(\Omega),
\]
\[
\int_\Omega \left(b(x,Du),D\phi \right)\,dx =\left\langle f,\phi
\right\rangle \quad \mbox{for every} \quad \phi \in H_0^{1,2}(\Omega ).
\]
Thus
\[
\int_\Omega \left(b^k\left(x,Du_{*}^k\right)-b^k(x,Du),D\phi
\right)\,dx= \int_\Omega
\left(b(x,Du)-b^k(x,Du),D\phi \right)\,dx,
\]
for every $\phi \in H_0^{1,2}(\Omega )$. Choose $\phi =u_{*}^k-u$ and take
the strict monotonicity of $b^k$, see Remark~1 in the next section,
into account on the left hand side and apply the Schwarz inequality and
H\"{o}lder inequality on the right hand side to obtain
\[
\alpha \int_\Omega \left| Du_{*}^k-Du\right| ^2\!\! dx\leq \left( \int_\Omega
\left| b(x,Du)-b^k(x,Du)\right| ^2\!\! dx\right) ^{\frac 12}\left( \int_\Omega
\left| Du_{*}^k-Du)\right| ^2\! dx\right) ^{\frac 12}\!.
\]
Hence, by Theorem~3
\begin{equation}
\left( \int_\Omega \left| Du_{*}^k-Du)\right| ^2\,dx\right) ^{\frac 12}\leq
\frac 1\alpha \left( \omega \left(\frac 1k\right)C\int_\Omega \left|
Du\right| ^2\,dx\right) ^{\frac 12}.  \label{step31}
\end{equation}
The right hand side tends to $0$ as $k\rightarrow \infty .$ We obtain (3) by
noting that $\left\| D\cdot \right\| _{L^2(\Omega ;R^n)}^2$ is an equivalent
norm on $H_0^{1,2}(\Omega )$.
\hfill {}$\square $

Next we prove that $\ds a\left(x,\frac x{\varepsilon
_h},Du_h\right)\rightarrow b(x,Du)$
weakly in $L^2(\Omega ;R^n)$. In fact if $g\in (L^2(\Omega ;R^n))^{*},$ then
\[
\ba{l}
\ds \lim_{h\rightarrow \infty }\left\langle g,a\left(x,\frac x{\varepsilon
_h},Du_h\right)-b(x,Du)\right\rangle =\lim_{k\rightarrow \infty
}\lim_{h\rightarrow \infty }\left\langle g,a\left(x,\frac x{\varepsilon
_h},Du_h\right)-b(x,Du)\right\rangle \\[5mm]
\ds \qquad =\lim_{k\rightarrow \infty }\lim_{h\rightarrow \infty }
\left\langle g,a\left(x,\frac x{\varepsilon _h},Du_h\right)-
a^k\left(x,\frac x{\varepsilon_h},Du_h^k\right)\right\rangle \\[5mm]
\ds \qquad +\lim_{k\rightarrow \infty }\lim_{h\rightarrow \infty }\left\langle
g,a^k\left(x,\frac x{\varepsilon
_h},Du_h^k\right)-b^k\left(x,Du_{*}^k\right) \right\rangle \\[5mm]
\ds \qquad +\lim_{k\rightarrow \infty }\lim_{h\rightarrow \infty }\left\langle
g,b^k(x,Du_{*}^k)-b(x,Du)\right\rangle \\[5mm]
\ds  \leq\lim_{k\rightarrow \infty }\lim_{h\rightarrow \infty }\left\|
g\right\| _{(L^2(\Omega ;R^n))^{*}}\left\| a\left(x,\frac x{\varepsilon
_h},Du_h\right)-a^k\left(x,\frac x{\varepsilon
_h},Du_h^k\right)\right\|_{L^2(\Omega ;R^n)}
\\[5mm]
\ds \qquad+\lim_{k\rightarrow \infty }\lim_{h\rightarrow \infty }\left\langle
g,a^k\left(x,\frac x{\varepsilon _h},Du_h^k\right)-
b^k\left(x,Du_{*}^k)\right)\right\rangle \\[5mm]
\ds \qquad+\lim_{k\rightarrow \infty }\left\| g\right\| _{(L^2(\Omega
;R^n))^{*}}\left\| b^k\left(x,Du_{*}^k\right)-b(x,Du)\right\|
_{L^2(\Omega ;R^n)}.
\ea
\]
It is suf\/f\/icient to prove that all three terms on the right hand side are
zero. We do this in three steps.

\textbf{Step 1.} Let us show that
\begin{equation}
\lim_{k\rightarrow \infty }\lim_{h\rightarrow \infty }\left\| a\left(x,\frac
x{\varepsilon _h},Du_h\right)-a^k\left(x,\frac x{\varepsilon
_h},Du_h^k\right) \right\|_{L^2(\Omega ;R^n)}=0.  \label{step12}
\end{equation}
By using elementary estimates we f\/ind that
\[
\ba{l}
\ds \int_\Omega \left| a^k\left(x,\frac x{\varepsilon
_h},Du_h^k\right) -a\left(x,\frac x{\varepsilon _h},Du_h\right)
\right| ^2\,dx \\[5mm]
\ds =\int_\Omega \left| a^k\!\left(x,\frac x{\varepsilon
_h},Du_h^k\right) -a^k\!\left(x,\frac x{\varepsilon _h},Du_h\right)
+a^k\!\left(x,\frac x{\varepsilon _h},Du_h\right)-a\!\left(x,\frac
x{\varepsilon _h},Du_h\right)\right| ^2\!\! dx \\[5mm]
\ds \leq 2\int_\Omega \left| a^k\left(x,\frac x{\varepsilon
_h},Du_h^k\right) -a^k\left(x,\frac x{\varepsilon _h},Du_h\right)
\right| ^2\,dx \\[5mm]
\ds +2\int_\Omega \left| a^k\left(x,\frac x{\varepsilon
_h},Du_h\right) -a\left(x,\frac x{\varepsilon _h},Du_h\right)\right| ^2\,dx.
\ea
\]
Hence, by applying the continuity conditions (\ref{acont}) and (\ref{cont3}%
), we obtain that
\[
\ba{l}
\ds \int_\Omega \left| a^k\left(x,\frac x{\varepsilon
_h},Du_h^k\right)- a\left(x,\frac x{\varepsilon
_h},Du_h\right)\right|^2\,dx \\[5mm]
\ds \qquad \leq 2\beta ^2\int_\Omega \left| Du_h^k-Du_h\right|
^2\,dx+ 2\omega \left(\frac 1k\right)\int_\Omega \left| Du_h\right| ^2\,dx.
\ea
\]
According to (\ref{bound}) and the fact that $(Du_h)$ is bounded in
$L^2(\Omega ;R^n)$ there exists a constant~$C$ independent of $h$
such that
\[
\int_\Omega \left| a^k\left(x,\frac x{\varepsilon _h},Du_h^k\right)-
a\left(x,\frac x{\varepsilon _h},Du_h\right)\right| ^2\,dx\leq
C\omega \left(\frac 1k\right).
\]
By the properties of $\omega $ it follows that
\begin{equation}
\left\| a\left(x,\frac x{\varepsilon _h},Du_h\right)-
a^k\left(x,\frac x{\varepsilon _h},Du_h^k\right)
\right\| _{L^2(\Omega ;R^n)}\rightarrow 0,  \label{step121}
\end{equation}
as $k\rightarrow \infty $ uniformly in $h.$ This implies that we may change
the order in the limit process in (\ref{step12}) and we obtain
(\ref{step12}) by taking (\ref{step121}) into account.
\hfill {}$\square $

\textbf{Step 2.} We observe that
\[
\lim_{k\rightarrow \infty }\lim_{h\rightarrow \infty }\left\langle
g,a^k\left(x,\frac x{\varepsilon _h},Du_h^k\right)-
b^k\left(x,Du_{*}^k\right)\right\rangle =0,
\]
as a direct consequence of Theorem~1. \hfill {}$\square $

\textbf{Step 3.} Let us show that
\begin{equation}
\lim_{k\rightarrow \infty }\left\| b^k\left(x,Du_{*}^k\right)-
b(x,Du)\right\| _{L^2(\Omega ;R^n)}=0.  \label{step13}
\end{equation}
We have that
\[
\ba{l}
\ds \int_\Omega \left| b^k\left(x,Du_{*}^k\right)-b(x,Du)\right|
^2\,dx \\[5mm]
\ds \qquad =\int_\Omega \left| b^k\left(x,Du_{*}^k\right)-
b^k(x,Du)+b^k(x,Du)-b(x,Du)\right|
^2\,dx \\[5mm]
\ds \qquad \leq 2\int_\Omega \left| b^k\left(x,Du_{*}^k\right)-
b^k(x,Du)\right|
^2\,dx+2\int_\Omega \left| b^k(x,Du)-b(x,Du)\right| ^2\,dx.
\ea
\]
By applying the continuity conditions in Remark~1 and Theorem~3 we see that
\[
\!\!\int_\Omega \left| b^k\left(x,Du_{*}^k\right)-b(x,Du)\right|
^2\,dx\leq 2 \frac{\beta ^4}{\alpha ^2}\int_\Omega \left|
Du_{*}^k-Du\right| ^2\,dx+2C\omega \left(\frac 1k\right)\int_\Omega
\left| Du\right| ^2\,dx.
\]
Now (\ref{step13}) follows by taking (\ref{step31}) into account.
\hfill {}$\square $

\section{Properties of the homogenized operator}

In this section we prove some properties of the homogenized operator. In
particular these properties implies the existence and uniqueness of the
solution of the homogenized problem.

\medskip

\noindent
{\bf Theorem 3.} \label{thbprop}
{\it Let $b$ be the homogenized operator defined in Theorem~2. Then

\begin{enumerate}
\item[(a)]  $b(\cdot ,\xi )$ satisfies the continuity condition
\[
\left| b(x_1,\xi )-b(x_2,\xi )\right| ^2\leq \omega (\left| x_1-x_2\right|
)C\left| \xi \right| ^2,
\]
where $\ds C=2\left( \frac \beta \alpha \right) ^2\left( 1+\left( \frac \beta
\alpha \right) ^2\right) $.

\item[(b)]  $b(x,\cdot )$ is strictly monotone, more precisely
\[
(b(x,\xi _1)-b(x,\xi _2),\xi _1-\xi _2)\geq \alpha \left| \ \xi _1-\xi
_2\right| ^2,
\quad \mbox{for every} \quad  \xi _1,\xi _2\in R^n.
\]

\item[(c)]  $b(x,\cdot )$ is Lipschitz continuous, more precisely
\[
\left| b(x,\xi _1)-b(x,\xi _2)\right| \leq \frac{\beta ^2}\alpha \left| \
\xi _1-\xi _2\right| ,
\quad \mbox{for every} \quad \xi _1,\xi _2\in R^n.
\]

\item[(d)]  $b(x,0)=0$ for $x\in \Omega .$
\end{enumerate}
}

\noindent
{\bf Proof.}
(a) By the def\/inition of $b$ we have
\[
\ba{l}
\ds \left| b(x_1,\xi )-b(x_2,\xi )\right| ^2 =\left|
\int_Ya\left(x_1, y,\xi +Dv^{\xi ,x_1}\right)\,dy-\int_Y
a\left(x_2,y,\xi +Dv^{\xi ,x_2}\right)\,dy\right| ^2 \\[5mm]
\ds \qquad \leq \Bigl(\int_Y\left| a\left(x_1,y,\xi +Dv^{\xi
,x_1}\right) -a\left(x_2,y,\xi +Dv^{\xi
,x_1}\right)\right| \,dy \\[5mm]
\ds \qquad +\int_Y\left| a\left(x_2,y,\xi +Dv^{\xi ,x_1}\right)-
a\left(x_2,y,\xi +Dv^{\xi ,x_2}\right)\right| \,dy\Bigr)^2 \\[5mm]
\ds \qquad \leq 2\left( \int_Y\left| a\left(x_1,y,\xi +Dv^{\xi
,x_1}\right)- a\left(x_2,y,\xi +Dv^{\xi ,x_1}\right)\right|
\,dy\right) \\[5mm]
\ds \qquad +2\left( \int_Y\left| a\left(x_2,y,\xi +Dv^{\xi ,x_1}\right)-
a\left(x_2,y,\xi +Dv^{\xi ,x_2}\right)\right| \,dy\right) ^2,
\ea
\]
where we in the last inequality used that $(a+b)^2\leq 2(a^2+b^2)$
for $a,b\geq 0$. By the Jensen inequality we get
\[
\ba{l}
\left| b(x_1,\xi )-b(x_2,\xi )\right| ^2 \leq 2\int_Y\left|
a\left(x_1,y,\xi +Dv^{\xi ,x_1}\right)-
a\left(x_2,y,\xi +Dv^{\xi ,x_1}\right)\right| ^2\,dy \\[5mm]
\ds \qquad +2\int_Y\left| a\left(x_2,y,\xi +Dv^{\xi ,x_1}\right)-
a\left(x_2,y,\xi +Dv^{\xi ,x_2}\right)\right| ^2\,dy.
\ea
\]

From the continuity conditions (\ref{cont3}) and (\ref{acont}) it
follows that
\begin{equation}
\ba{l}
\left| b(x_1,\xi )-b(x_2,\xi )\right| ^2\\[5mm]
\ds \qquad \leq 2\omega (\left| x_1-x_2\right|
)\int_Y\left| \xi +Dv^{\xi ,x_1}\right| ^2\,dy+2\beta ^2\int_Y\left| Dv^{\xi
,x_1}-Dv^{\xi ,x_2}\right| ^2\,dy.
\ea
\label{bcontequ1}
\end{equation}
We will now study the two terms in (\ref{bcontequ1}) separately.

First term: According to (\ref{amon}), (\ref{cellproblem}), (\ref{acont}) we
have
\[
\ba{l}
\ds \alpha \int_Y\left| \xi _1+Dv^{\xi _1,x}(y)-\xi _2-Dv^{\xi
_2,x}(y)\right| ^2\,dy \\[5mm]
\ds \leq \int_Y\left(a\left(x,y,\xi _1+Dv^{\xi _1,x}\right)-
a\left(x,y,\xi _2+Dv^{\xi _2,x}\right),\xi
_1+Dv^{\xi _1,x}-\xi _2-Dv^{\xi _2,x}\right)\,dy \\[5mm]
\ds =\int_Y\left(a\left(x,y,\xi _1+Dv^{\xi _1,x}\right)-a\left(x,y,\xi
_2+Dv^{\xi _2,x}\right),\xi _1-\xi  _2\right)\,dy \\[5mm]
\ds \leq \int_Y\left| a\left(x,y,\xi _1+Dv^{\xi _1,x}\right)-
a\left(x,y,\xi _2+Dv^{\xi _2,x}\right)\right| \left| \xi _1-\xi
_2\right| \,dy \\[5mm]
\ds \leq \beta \int_Y\left| \xi _1+Dv^{\xi _1,x}-\xi _2-Dv^{\xi _2,x})\right|
\left| \xi _1-\xi _2\right| \,dy \\[5mm]
\ds \leq \beta \left( \int_Y\left| \xi _1+Dv^{\xi _1,t}-\xi _2-Dv^{\xi
_2,t})\right| ^2\,dy\right) ^{\frac 12}\left( \int_Y\left| \xi _1-\xi
_2\right| ^2\,dy\right) ^{\frac 12}.
\ea
\]
This can be written as
\begin{equation}
\int_Y\left| \xi _1+Dv^{\xi _1,t}(y)-\xi _2-Dv^{\xi _2,t}(y)\right|
^2\,dy\leq \left( \frac \beta \alpha \right) ^2\left| \xi _1-\xi _2\right|
^2.  \label{corrcont}
\end{equation}

Second term: By def\/inition we have that
\[
\ba{l}
\ds \int_Y\left(a\left(x_1,y,\xi +Dv^{\xi ,x_1}(y)\right),D\phi
\right)\,dy =0
\quad \mbox{for every} \quad \phi
\in H_{\square }^{1,2}(\Omega ), \\[5mm]
\ds \int_Y\left(a\left(x_2,y,\xi +Dv^{\xi ,x_2}(y)\right),D\phi
\right)\,dy =0
\quad \mbox{for every} \quad \phi
\in H_{\square }^{1,2}(\Omega ).
\ea
\]
This implies that
\[
\ba{l}
\ds \int_Y\left(a\left(x_1,y,\xi +Dv^{\xi ,x_1}\right)-
a\left(x_1,y,\xi +Dv^{\xi ,x_2}\right),D\phi \right)\,dy
\\[5mm]
\ds \qquad =\int_Y\left(a\left(x_2,y,\xi +Dv^{\xi ,x_2}\right)-
a\left(x_1,y,\xi +Dv^{\xi ,x_2}\right),D\phi \right)\,dy,
\ea
\]
for every $\phi \in H_{\square }^{1,2}(\Omega )$. In particular for $\phi
=Dv^{\xi ,x_1}-Dv^{\xi ,x_2}$ we have
\[
\ba{l}
\ds \int_Y\left(a\left(x_1,y,\xi +Dv^{\xi ,x_1}\right)-
a\left(x_1,y,\xi +Dv^{\xi ,x_2}\right),\xi
+Dv^{\xi ,x_1}-\left(\xi +Dv^{\xi ,x_2}\right)\right)\,dy \\[5mm]
\ds \qquad =\int_Y\left(a\left(x_2,y,\xi +Dv^{\xi ,x_2}\right)-
a\left(x_1,y,\xi +Dv^{\xi ,x_2}\right),Dv^{\xi
,x_1}-Dv^{\xi ,x_2}\right)\,dy.
\ea
\]
By applying (\ref{amon}) on the left hand side and the Schwarz and
H\"{o}lder inequalities on the right hand side it follows that
\[
\ba{l}
\ds \alpha \left( \int_Y\left| Dv^{\xi ,x_1}-Dv^{\xi ,x_2}\right| ^2\,dy\right)
^{\frac 12}\\[5mm]
\ds \qquad \leq \left( \int_Y\left| a\left(x_2,y,\xi +Dv^{\xi
,x_2}\right) -a\left(x_1,y,\xi
+Dv^{\xi ,x_2}\right)\right| ^2\,dy\right) ^{\frac 12}.
\ea
\]
The continuity condition (\ref{cont3}) and (\ref{bcontequ2}) imply that
\be
\ba{l}
\ds \alpha \left( \int_Y\left| Dv^{\xi ,x_1}-Dv^{\xi ,x_2}\right|
^2\,dy \right) ^{\frac 12} \\[5mm]
\ds \qquad \leq \left( \omega (\left| x_1-x_2\right| )\int_Y\left| \xi
+Dv^{\xi ,x_2}\right| ^2\,dy\right) ^{\frac 12}
\leq \omega ^{\frac 12}(\left| x_1-x_2\right| )\frac \beta \alpha \left|
\xi \right| .
\ea  \label{bcontequ2}
\ee
The result follows by taking (\ref{bcontequ1}), (\ref{corrcont}) and (\ref
{bcontequ2}) into account.

(b) Let $\xi _j\in R^n$, $j=1,2$ and def\/ine \thinspace for a.e. $y\in R^n,$
\[
w_h^{\xi _j,x}(y)=(\xi _j,y)+\varepsilon _hv^{\xi _j,x}\left(\frac
y{\varepsilon _h}\right).
\]
By the periodicity of $v^{\xi _j,x}$ we have that
\be
w_h^{\xi _j,x} \rightarrow (\xi _j,y)
\quad \mbox{weakly in} \quad H^{1,2}(\Omega ),
\label{conv1}
\ee
\be
Dw_h^{\xi _j,x} \rightarrow \xi _j \quad \mbox{weakly in}
\quad L^2(\Omega ,R^n),
\label{conv2}
\ee
\be
a\left(x,\frac y{\varepsilon _h},Dw_h^{\xi _j,x}(y)\right)
\rightarrow b(x,\xi _j) \quad \mbox{weakly in} \quad L^2(\Omega
,R^n).   \label{conv3}
\ee
The monotonicity condition (\ref{amon}) on $a$ implies that
\[
\ba{l}
\ds \int_\Omega \left(a\left(x,\frac y{\varepsilon _h},Dw_h^{\xi
_1,x}\right) -a\left(x,\frac y{\varepsilon _h},Dw_h^{\xi
_2,x}\right), Dw_h^{\xi _1,x}-Dw_h^{\xi _2,x}\right)\phi
(y)\,dy \\[5mm]
\ds \qquad \geq \alpha \int_\Omega \left| Dw_h^{\xi _1,x}-Dw_h^{\xi
_2,x}\right| ^2\phi (y)\,dy
\ea
\]
for every $\phi \in C_0^\infty (\Omega )$ such that $\phi \geq 0$. We
apply $\lim \inf\limits_{h\rightarrow \infty }$ on both sides of this
inequality and obtain
\[
\int_\Omega (b(x,\xi _1)-b(x,\xi _2),\xi _1-\xi _2)\phi (y)\,dy\geq \alpha
\int_\Omega \left| \xi _1-\xi _2\right| ^2\phi (y)\,dy,
\]
for every $\phi \in C_0^\infty (\Omega )$, where we on the left hand side
have used (\ref{conv2}), (\ref{conv3}), (\ref{divwhk}) and Lemma~1
and on the right hand side we have used
(\ref{conv2}) and the fact that for a weakly convergent sequence
$(x_n)$ converging to $x$ we have that $\left\| x\right\| \leq \lim
\inf\limits_{n\rightarrow \infty
}\left\| x_n\right\| $. This implies that
\[
(b_i(\xi _1)-b_i(\xi _2),\xi _1-\xi _2)\geq \alpha \left| \ \xi _1-\xi
_2\right| ^2.
\]
Since $\xi _1,\xi _2$ were chosen arbitrary (b) is proved.

(c) Fix $\xi _1,\xi _2\in R^n$. According to (\ref{acont}), Jensen\'{}s
inequality and (\ref{amon}) we have that
\[
\ba{l}
\hspace*{-9.7pt}\ds \left| b(x,\xi _1)-b(x,\xi _2)\right| ^2
 =\left| \int_Ya\left(x,y,\xi
_1+Dv^{\xi _1,x}(y)\right)\! dy-\!\int_Ya\left(x,y,\xi _2+Dv^{\xi
_2,x}(y)\right) \! dy\right| ^2
\\[5mm]
\ds \qquad \leq \left( \int_Y\left| a\left(x,y,\xi _1+Dv^{\xi
_1,x}(y)\right) -a\left(x,y,\xi _2+Dv^{\xi _2,x}(y)\right)\right| \,
dy\right) ^2 \\[5mm]
\ds \qquad \leq \left( \int_Y\beta \left| \xi _1+Dv^{\xi _1,x}(y)-
\xi _2-Dv^{\xi_2,x}(y)\right| \,dy\right) ^2 \\[5mm]
\ds \qquad \leq \beta ^2\int_Y\left| \xi _1+Dv^{\xi _1,x}(y)-\xi _2-Dv^{\xi
_2,x}(y)\right| ^2\,dy \\[5mm]
\ds \qquad \leq \frac{\beta ^2}\alpha \int_Y\left(a\left(x,y,\xi _1+Dv^{\xi
_1,x}(y)\right)-a\left(x,y,\xi  _2+Dv^{\xi _2,x}(y)\right), \right.\\[5mm]
\left.\ds \qquad \left(\xi _1+Dv^{\xi _1,x}(y)\right)-\left(\xi _2+Dv^{\xi
_2,x} (y)\right)\right)\,dy.
\ea
\]
Moreover, from the cell problem (\ref{cellproblem}) it follows that
\[
\int_Y\left(a\left(x,y,\xi _k+Dv^{\xi _k,x}(y)\right),Dv^{\xi
_l,x}(y)\right) \,dx=0
\quad \mbox{for} \quad
k,l=1,2,
\]
and we conclude
\[
\ba{l}
\ds \left| b(x,\xi _1)-b(x,\xi _2)\right| ^2\\[5mm]
\ds \qquad \leq \frac{\beta ^2}\alpha (b(x,\xi
_1)-b(x,\xi _2),\xi _1-\xi _2)\leq \frac{\beta ^2}\alpha \left| b(x,\xi
_1)-b(x,\xi _2)\right| \left| \xi _1-\xi _2\right|
\ea
\]
which implies (c).

(d) Since $a(x,y,0)=0$ we have that the solution of the cell-problem (\ref
{cellproblem}) corresponding to $\xi =0$ is $v^{0,x}=0$. This implies that
\[
b(x,0)=\int_Ya(x,y,0)\,dy=0. \qquad \mbox{\hspace*{8.0cm}  \hfill $\square $}
\]

\noindent
{\bf Remark 1.} \label{rem}
By similar arguments it follows that (b), (c) and (d) holds, up
to boundaries, for the homogenized operator $b$ in Theorem~1.

\section{Some corrector results}

We have proved in both Theorem~1 and Theorem~2 that we
for the corresponding solutions have that $u_u-u$ converges to $0$ weakly in
$H_0^{1,2}(\Omega )$. By the Rellich imbedding theorem we have that $u_h-u$
converges to $0$ in $L^2(\Omega )$. In general we do not have strong
convergence of $Du_h-Du$ to $0$ in $L^2(\Omega ,R^n)$. However, we will
prove that it is possible to express $Du_h$ in terms of $Du$, up to a rest
which converges to $0$ in $L^2(\Omega ,R^n)$.

\medskip

\noindent
{\bf Theorem 4.} \label{thtrancor}
{\it Let $u$ and $u_h$ be defined as in Theorem~1
and let $P_h$ be given by (\ref{pdef}). Then
\[
Du_h-\sum_{i=1}^N\chi _{\Omega _i}(x)P_h(x,M_hDu,x_i)\rightarrow 0
\quad \mbox{in} \quad L^2(\Omega ,R^n).
\]
}

\noindent
{\bf Proof.}
In \cite{dalmaso1992} the case $N=1$ was considered. By using these ideas
and make the necessary adjustments the proof follows. For the details the
reader is referred to \cite{wall1998c}.  \hfill {}$\square $

\medskip

\noindent
{\bf Theorem 5.} \label{thcontcor}
{\it Let $u$ and $u_h$ be defined as in Theorem~2.
 Moreover, let $P_h$ be given by (\ref{pdef}) and $\gamma _h$ by
(\ref{gamma}). Then
\[
Du_h-P_h(x,M_hDu,\gamma _h)\rightarrow 0
\quad \mbox{in} \quad L^2(\Omega ,R^n)
\]
}

\noindent
{\bf Proof.} We have that
\be \label{equ1thconcor}
\ba{l}
\ds \left\| Du_h-M_hDu-Dv^{M_hDu,\gamma _h}
\left(\frac x{\varepsilon _h}\right)\right\|
_{L^2(\Omega ,R^n)}  \leq \left\| Du_h-Du_h^k\right\| _{L^2(\Omega ,R^n)}
\\[5mm]
\ds +\left\|
Du_h^k-M_hDu_{*}^k-\sum_{i\in I_k}\chi _{\Omega
_i^k}(x)Dv^{M_hDu_{*}^k,x_i^k}\left(\frac x{\varepsilon _h}\right)
\right\| _{L^2(\Omega ,R^n)} \\[7mm]
\ds +\left\| M_hDu_{*}^k+\sum_{i\in I_k}\chi _{\Omega
_i^k}Dv^{M_hDu_{*}^k,x_i^k}\left(\frac x{\varepsilon _h}\right)-
M_hDu-Dv^{M_hDu,\gamma
_h}\left(\frac x{\varepsilon _h}\right)\right\| _{L^2(\Omega ,R^n)}.
\ea\hspace{-6pt}
\ee
As in the proof of Theorem~2 we have that
\[
\lim_{k\rightarrow \infty }\lim_{h\rightarrow \infty }\left\|
Du_h-Du_h^k\right\| _{L^2(\Omega ,R^n)}=0
\]
and by Theorem~4 it yields that
\[
\lim_{k\rightarrow \infty }\lim_{h\rightarrow \infty }\left\|
Du_h^k-M_hDu_{*}^k-\sum_{i\in I_k}\chi _{\Omega
_i^k}(x)Dv^{M_hDu_{*}^k,x_i^k}\left(\frac x{\varepsilon _h}\right)
\right\| _{L^2(\Omega,R^n)}=0.
\]
This means that the theorem would be proved if we show that
$\lim\limits_{k\rightarrow \infty }\lim\limits_{h\rightarrow \infty }$
acting on the last
term in (\ref{equ1thconcor}) is equal to $0$. In order to prove this fact we
make the following elementary estimations:
\be
\ba{l}
\ds \left\| M_hDu_{*}^k+\sum_{i\in I_k}\chi _{\Omega
_i^k}(x)Dv^{M_hDu_{*}^k,x_i^k}\left(\frac x{\varepsilon
_h}\right)-M_hDu-Dv^{M_hDu,\gamma _h}\left(\frac x{\varepsilon _h}
\right)\right\| _{L^2(\Omega ,R^n)}^2 \\[6mm]
\ds =\sum_{i\in I_k}\int_{\Omega _i^k}\left|
M_hDu_{*}^k+Dv^{M_hDu_{*}^k,x_i^k}\left(\frac x{\varepsilon
_h}\right)-M_hDu-Dv^{M_hDu,\gamma _h}\left(\frac x{\varepsilon _h}
\right)\right| ^2\,dx  \\[6mm]
\ds \leq \sum_{i\in I_k}\int_{\Omega _i^k}\left( \left|
M_hDu_{*}^k+Dv^{M_hDu_{*}^k,x_i^k}\left(\frac x{\varepsilon
_h}\right)-M_hDu-Dv^{M_hDu,x_i^k}\left(\frac x{\varepsilon _h}\right)
\right| \right.  \\[6mm]
\ds \qquad \left. +\left| Dv^{M_hDu,x_i^k}\left(\frac x{\varepsilon _h}
\right)-Dv^{M_hDu,\gamma _h}\left(\frac x{\varepsilon _h}\right)
\right| \right) ^2\,dx  \\[6mm]
\ds \leq \sum_{i\in I_k}2\int_{\Omega _i^k}\left|
M_hDu_{*}^k+Dv^{M_hDu_{*}^k,x_i^k}\left(\frac x{\varepsilon
_h}\right)-M_hDu-Dv^{M_hDu,x_i^k}\left(\frac x{\varepsilon _h}\right)
\right| ^2\,dx  \\[6mm]
\ds \qquad +\sum_{i\in I_k}2\int_{\Omega _i^k}\left| Dv^{M_hDu,x_i^k}
\left(\frac x{\varepsilon _h}\right)-
Dv^{M_hDu,\gamma _h}\left(\frac x{\varepsilon _h}\right)\right|^2\,dx.
\ea \label{equ2thconcor}\hspace{-5.5pt}
\ee
We shall now study the two terms on the right hand side of
(\ref{equ2thconcor}) separately but f\/irst we def\/ine,
\[
\xi _{h,*}^{j,k}=\frac 1{\left| Y_h^j\right| }\int_{Y_h^j}Du_{*}^k\,dx.
\]
By using (\ref{corrcont}) and a change of variables we f\/ind that we for the
f\/irst term in (\ref{equ2thconcor}) have the following estimate
\be \label{equ3thconcor} \\
\ba{l}
\ds \int_{\Omega _i^k}\left| M_hDu_{*}^k+Dv^{M_hDu_{*}^k,x_i^k}
\left(\frac x{\varepsilon _h}\right)
-M_hDu-Dv^{M_hDu,x_i^k}\left(\frac x{\varepsilon _h}\right)\right|
^2\,dx   \\[6mm]
\ds \leq \sum_{j\in J_h^{i,k}}\int_{Y_h^j}\left| \xi _{h,*}^{j,k}+
Dv^{\xi _{h,*}^{j,k},x_i^k}\left(\frac x{\varepsilon _h}\right)
-\xi _h^j-Dv^{\xi _h^j,x_i^k}\left(\frac x{\varepsilon _h}\right)
\right| ^2\,dx  \\[6mm]
\ds \qquad+\sum_{j\in B_h^{i,k}}\int_{Y_h^j}\left| \xi _{h,*}^{j,k}+Dv^{\xi
_{h,*}^{j,k},x_i^k}\left(\frac x{\varepsilon _h}\right)
-\xi _h^j-Dv^{\xi _h^j,x_i^k}\left(\frac x{\varepsilon _h}\right)
\right| ^2\,dx  \\[6mm]
\ds \leq \sum_{j\in J_h^{i,k}}\left( \frac \beta \alpha \right)^2\left| \xi
_{h,*}^{j,k}-\xi _h^j\right| ^2\left| Y_h^j\right| +\sum_{j\in
B_h^{i,k}}\left( \frac \beta \alpha \right) ^2\left| \xi _{h,*}^{j,k}-\xi
_h^j\right| ^2\left| Y_h^j\right|  \\[6mm]
\ds =\left( \frac \beta \alpha \right) ^2\int_{\Omega _i^{k,h}}\left|
M_hDu_{*}^k-M_hDu\right| ^2\,dx+\left( \frac \beta \alpha \right)
^2\int_{F_i^{k,h}}\left| M_hDu_{*}^k-M_hDu\right| ^2\,dx  \\[6mm]
\ds \leq \left( \frac \beta \alpha \right) ^2\int_{\Omega _i^k}\left|
Du_{*}^k-Du\right| ^2\,dx+\left( \frac \beta \alpha \right)
^2\int_{F_i^{k,h}}\left| Du_{*}^k-Du\right| ^2\,dx,
\ea
\ee
where we used Jensen\'{}s inequality in the last step. Moreover, by
using (\ref{bcontequ2}) and Jensen\'{}s inequality, we obtain that
\be \label{equ4thconcor}
\ba{l}
\ds \int_{\Omega _i^k}\left| Dv^{M_hDu,x_i^k}\left(\frac x{\varepsilon
_h}\right)-Dv^{M_hDu,\gamma _h}
\left(\frac x{\varepsilon _h}\right)\right| ^2\,dx  \\[6mm]
\qquad \ds \leq \sum_{j\in J_h^{i,k}}\int_{Y_h^j}\left| Dv^{\xi _h^j,x_i^k}
\left(\frac x{\varepsilon _h}\right)-
Dv^{\xi _h^j,x_h^j}\left(\frac x{\varepsilon _h}\right)\right| ^2\,dx
 \\[6mm]
\qquad \ds +\sum_{j\in B_h^{i,k}}\int_{Y_h^j}\left| Dv^{\xi _h^j,x_i^k}
\left(\frac x{\varepsilon _h}\right)-
Dv^{\xi _h^j,x_h^j}\left(\frac x{\varepsilon _h}\right)\right| ^2\,dx
 \\[6mm]
\qquad \ds \leq \sum_{j\in J_h^{i,k}}\omega \left(\frac 1k\right)
\frac{\beta ^2}{\alpha ^4}
\left| \xi _h^j\right| ^2\left| Y_h^j\right| +\sum_{j\in B_h^{i,k}}\omega
\left(\frac 1k+n\varepsilon _h\right)\frac{\beta ^2}{\alpha ^4}
\left| \xi _h^j\right| ^2\left| Y_h^j\right| \\[6mm]
\qquad \ds =\omega \left(\frac 1k\right)\frac{\beta ^2}{\alpha ^4}
\int_{\Omega _i^{k,h}}\left|
M_hDu\right| ^2\,dx+\omega \left(\frac 1k+n\varepsilon _h\right)
\frac{\beta ^2}{\alpha
^4}\int_{F_i^{k,h}}\left| M_hDu\right| ^2\,dx  \\[6mm]
\qquad \ds \leq \omega \left(\frac 1k\right)
\frac{\beta ^2}{\alpha ^4}\int_{\Omega _i^k}\left|
Du\right| ^2\,dx+\omega \left(\frac 1k+n\varepsilon _h\right)
\frac{\beta ^2}{\alpha ^4}
\int_{F_i^{k,h}}\left| Du\right| ^2\,dx.
\ea
\ee
By combining (\ref{equ2thconcor}), (\ref{equ3thconcor}) and
(\ref{equ4thconcor}) we have that
\[
\ba{l}
\ds \left\| M_hDu_{*}^k+\sum_{i\in I_k}\chi _{\Omega
_i^k}Dv^{M_hDu_{*}^k,x_i^k}\left(\frac x{\varepsilon _h}\right)
-M_hDu-Dv^{M_hDu,\gamma
_h}\left(\frac x{\varepsilon _h}\right)\right\| _{L^2(\Omega ,R^n)}^2 \\[6mm]
\qquad \ds \leq 2\left( \frac \beta \alpha \right) ^2\left( \int_\Omega \left|
Du_{*}^k-Du\right| ^2\,dx+\sum_{i\in I_k}\int_{F_i^{k,h}}\left|
Du_{*}^k-Du\right| ^2\,dx\right)  \\[6mm]
\qquad \ds +2\omega \left(\frac 1k+n\varepsilon _h\right)
\frac{\beta ^2}{\alpha ^4}\left(
\int_\Omega \left| Du\right| ^2\,dx+\sum_{i\in I_k}\int_{F_i^{k,h}}\left|
Du\right| ^2\,dx\right) .
\ea
\]
Moreover, by noting that $\left| F_i^{k,h}\right| \rightarrow 0$ as
$h\rightarrow \infty $ and taking (\ref{step3}) into account, we obtain that
$$
\ds \lim_{k\rightarrow \infty }\lim_{h\rightarrow \infty }\left\|
M_hDu_{*}^k+\sum_{i\in I_k}\chi _{\Omega _i^k}Dv^{M_hDu_{*}^k,x_i^k}
\left(\frac x{\varepsilon _h}\right)
 -M_hDu-Dv^{M_hDu,\gamma _h}
\left(\frac x{\varepsilon _h}\right)\right\| _{L^2(\Omega ,R^n)}\!\!\!\!=0
$$
and the theorem is proved.
\hfill {}$\square $

\subsection*{Acknowledgment}
The author thank the referee for some advises which
have improved the f\/inal version of this paper.

\label{wall-lp}

\end{document}